\documentclass{article}[12pt]
\usepackage[utf8]{inputenc} 
\usepackage[T1]{fontenc}
\usepackage{graphicx}
\usepackage{amsmath}
\usepackage{amssymb}
\usepackage{amsthm}
\usepackage{dsfont}
\usepackage{tikz}
\usepackage{authblk}
\usepackage{a4wide}

\newtheorem{theo}{Theorem}
\newtheorem{prop}{Proposition}
\newtheorem{defi}{Definition}

\newtheorem{cor}{Corollary}

\title{Generalizations of Efron's theorem}
\author{Yannis Oudghiri \footnote{yannis.oudghiri@univ-amu.fr}}
\affil{Aix Marseille Univ, CNRS, Centrale Marseille, I2M, Marseille, France}

\begin{document}

\maketitle

\begin{abstract}
In this article, we prove two new versions of a theorem proven by Efron in [Efr65]. Efron's theorem says that if a function $\phi : \mathbb{R}^2 \rightarrow \mathbb{R}$ is non-decreasing in each argument then we have that the function $s \mapsto \mathbb{E}[\phi(X,Y)|X+Y=s]$ is non-decreasing. We name restricted Efron's theorem a version of Efron's theorem where $\phi : \mathbb{R} \rightarrow \mathbb{R}$ only depends on one variable. \\
$PF_n$ is the class of functions such as 
$\forall a_1 \leq ... \leq a_n, b_1 \leq ... \leq b_n, \det(f(a_i-b_j))_{1 \leq i,j \leq n} \geq 0.$ The first version generalizes the restricted Efron’s theorem for random variables in the $PF_n$ class. The second one considers the non-restricted Efron’s theorem with a stronger monotonicity assumption. In the last part, we give a more general result of the second generalization of Efron's theorem.
\end{abstract}

\section{Introduction and notations}

As this article deals with log-concave random variables and probability density functions, we first recall basic definitions and properties:

\begin{defi}
A density $f : \mathbb{R} \rightarrow \mathbb{R}_+$ is said to be log-concave if it satisfies the following inequality:
\begin{equation*}
\forall x,y \in \mathbb{R}, \forall \theta \in [0,1], f(\theta x + (1 - \theta) y) \geq f(x)^{\theta}  f(y)^{1 - \theta}.
\end{equation*}
A real-valued random variable $X$ is said to be log-concave if it admits a density $f$ which is log-concave.
\end{defi}

 Basic properties of log-concave distributions are surveyed in a recent paper by Saumard and Wellner, see \cite{review}.

\medskip

Log-concavity for integer-valued probability distributions, i.e. for families ${(f(k))}_{k \in \mathbb{Z}}$ of non-negative numbers with $\underset{k \in \mathbb{Z}}{\sum} f(k) = 1$, can be defined as follows:
\begin{defi}
    We say that ${(f(k))}_{k \in \mathbb{Z}}$ is log-concave if we have the following inequality:
    \begin{equation*}
        \forall k \in \mathbb{Z} \ , \  f(k)^2 \geq f(k-1) f(k+1).
    \end{equation*}
\end{defi}

\medskip

\begin{defi}
Given an integer $n \ge 2$, function $f : \mathbb{R} \rightarrow \mathbb{R}$ is said to be a P\'olya frequency function in the class $PF_n$ if for every $n$-uples $a_1 \leq \dots \leq a_n$ and $b_1 \leq \dots \leq b_n$ we have:
\begin{equation*}
\det {(f(a_i - b_j))}_{1 \leq i,j \leq n} \geq 0.
\end{equation*}
\end{defi}

For example, the density $x \mapsto \dfrac{1}{\sqrt{2 \pi}} \exp \left( - \dfrac{x^2}{2}\right)$ is in the class $PF_n$ for any $n \geq 2$\cite{schoenberg}. \\
We can find a characterization of $PF_3$ class and some examples in \cite{PF3}.

The notion of $PF_n$ class can be seen as a natural generalization of log-concavity due to the following:

\begin{prop}
\label{prop1}
The class $PF_2$ is the class of log-concave functions.
\end{prop}
 The seminal paper \cite{schoenberg} by Schoenberg surveys basic properties of Pólya frequency functions.

\medskip

In \cite{efron}, Efron studied a non-trivial monotonicity property for log-concave distributions:
\begin{defi}
    Let $X,Y$ be two real-valued random variables.
    \begin{itemize}
        \item The pair $(X,Y)$ is said to satisfy the restricted monotonicity property if for any non-decreasing function $\phi : \mathbb{R} \rightarrow \mathbb{R}$, the functions ${\Phi}_X,{\Phi}_Y: \mathbb{R} \rightarrow \mathbb{R}$, defined by 
        \begin{equation*}
            {\Phi}_X(s) := \mathbb{E}[\phi(X)|X+Y = s] \ , \ {\Phi}_Y(s) := \mathbb{E}[\phi(Y)|X+Y = s].
        \end{equation*}
        are also non-decreasing.
        \item The pair (X,Y) is said to satisfy the strong monotonicity property if for any function $\phi : \mathbb{R}^2 \rightarrow \mathbb{R}$ which is non-decreasing in each variable, the function $\Phi : \mathbb{R} \rightarrow \mathbb{R}$, defined by 
        \begin{equation*}
            \Phi(s) := \mathbb{E}[\phi(X,Y)|X+Y =s].
        \end{equation*}
        is also non-decreasing.
    \end{itemize}
\end{defi}

Efron's theorem is stated as follows:
\begin{theo}\label{thm_efron}
    Let $(X,Y)$ be a pair of real-valued random variables.
    \begin{itemize}
        \item The pair $(X,Y)$ satisfies the strong monotonicity property if and only if satisfies the restricted monotonicity property.
        \item If $X$ and $Y$ are log-concave, then they satisfy the restricted and strong monotonicity properties.
    \end{itemize}
\end{theo}

By using Efron's theorem one can prove the stability of log-concavity by convolution. In the article \cite{Wellner}, there is also a proof of the stability of  strong log-concavity by convolution. 

\medskip

A recent work on Efron's monotonicity property is the article \cite{efron_2} by Saumard and Wellner. In this paper, the authors describe conditions on a pair $(X,Y)$, other than log-concavity, for the monotonicity property to hold. 

\medskip

The purpose of this paper is to prove two others generalizations of Efron's theorem. These generalizations are actually suggested by the title of Efron's original paper, \textit{"Increasing properties of P\'olya frequency functions"}, which is a bit misleading since the paper only deals with log-concave distributions.

\medskip

In Section 2, we generalize a restricted version of Efron's Theorem, see Theorem \ref{theorem1}: a generalized restricted monotonicity property holds for pairs of $PF_n$ variables. A consequence of this theorem is another proof of the stability of the class $PF_n$ by convolution. The proof of Theorem \ref{theorem1} is very simple and only uses elementary properties of determinants. \\

In section 3, we reinforce the strong version of Efron's Theorem: we prove that for a stronger monotonicity assumption on the function $\phi$, we obtain a more powerful result. We can see Theorem \ref{thm_efron} as a particular case of Theorem \ref{theorem2}. \\

In section 4, we give a more general result of the section  3, where Theorem \ref{theorem2} and Theorem \ref{theorem2_2} are a particular case. We apply in particular this result to the function $\Phi : s \mapsto \mathbb{E}[XY|X+Y=s]$.

\section{A restricted Efron's theorem for $PF_n$ variables}

In this paragraph, we prove the following generalization of restricted Efron's theorem.

\begin{defi}
A $n$-uple of functions ${\phi}_1, \ldots ,{\phi}_n$ is in the class $GM_n$ if for every $x_1 \leq \ldots \leq x_n$, we have :
\begin{equation*}
\det {({\phi}_i(x_j))}_{1 \leq i,j \leq n} \geq 0.
\end{equation*}
\end{defi}

We can remark that, for a pair $(\phi_1,\phi_2)$ of positive functions, the condition $({\phi}_1,{\phi}_2) \in GM_2$ is equivalent to say that the function $\dfrac{{\phi}_2}{{\phi}_1}$ is non-decreasing. 

\begin{theo}
\label{theorem1}
We fix some $n \geq 2$. Let $X,Y$ be independent random variables with $PF_n$ densities $f,g$. Let $({\phi}_1, \dots ,{\phi}_n)$ be a n-uple of functions in the class $GM_n$. For $1 \leq k \leq n$ we define the function ${\Phi}_k : \mathbb{R} \rightarrow \mathbb{R}$ by
\begin{equation*}
    {\Phi}_k(s) := \mathbb{E}[{\phi}_k(X) | X+Y = s].
\end{equation*}
We assume that ${\Phi}_k$ is well-defined. Then the n-uple $({\Phi}_1, \ldots ,{\Phi}_n)$ is also in the class $GM_n$.
\end{theo}

Restricted Efron's theorem is the particular case where $n=2$, ${\phi}_1 = 1$ and ${\phi}_2$ is a non-decreasing function. In order to prove Theorem \ref{theorem1}, we will use two elementary propositions. \\
The first one is a well-known identity which can be traced back to Andreev, see \cite{Andreief} and \cite{forrester}.

\begin{prop}
\label{proposition1}
Let ${(f_{i,j}(x))}_{1 \leq i,j \leq n}$ be a $n \times n$ matrix of integrable functions. Then
\begin{equation*}
\det \left( \int f_{i,j}(x) dx \right) = \int_{\mathbb{R}^n} \det(f_{i,j}(x_i)) dx_1 \dots dx_n = \int_{\mathbb{R}^n} \det(f_{i,j}(x_j)) dx_1 \dots dx_n    .
\end{equation*}
\end{prop}

The second elementary proposition is proven by noticing that the set of $(x_1, \dots ,x_n) \in \mathbb{R}^n$ such that $x_i = x_j$ for some $i \neq j$ has its Lebesgue measure equal to zero. We denote by $\mathfrak{S}_n$ the set of permutations on $\{ 1, \dots , n\}$. 

\begin{prop}
\label{proposition2}
Let $h : \mathbb{R}^n \rightarrow \mathbb{R}$ be an integrable function. We have:
\begin{equation*}
    \int_{\mathbb{R}^n} h(x_1, \dots ,x_n) dx_1 \dots dx_n = \int_{x_1 < \ldots < x_n} \underset{\sigma \in \mathfrak{S}_n}{\sum} h(x_{\sigma(1)}, \ldots ,x_{\sigma(n)}) dx_1 \ldots dx_n.
\end{equation*}
\end{prop}

\begin{proof}[Proof of Theorem \ref{theorem1}]
Let us fix an ordered n-uple $s_1 < \dots < s_n$. For $1 \leq i,j \leq n$, we have:
\begin{equation*}
    {\Phi}_i(s_j) = \mathbb{E}[{\phi}_i(X)|X+Y=s_j] = \dfrac{\int_{\mathbb{R}} f(x_i) {\phi}_i(x_i) g(s_j- x_i ) dx_i}{\int_{\mathbb{R}} f(x) g(s_j - x) dx}.
\end{equation*}
We thus have 
\begin{equation}
\label{equation1}
    \det({\Phi}_i(s_j)) = \dfrac{1}{\overset{n}{\underset{j=1}{\prod}} \int_{\mathbb{R}} f(x) g(s_j - x)} \det \left( \int_{\mathbb{R}} f(x_i) {\phi}_i(x_i) g(s_j -x_i) dx_i \right).
\end{equation}
Theorem \ref{theorem1} will thus be proven if we show that the determinant of the right-hand of \eqref{equation1} is non-negative. We will denote this determinant by $D$. \\
By Proposition \ref{proposition1}, we have:
\begin{equation*}
    D = \int_{\mathbb{R}^n} \left( \overset{n}{\underset{i=1}{\prod}} f(x_i) \right) \left( \overset{n}{\underset{i=1}{\prod}} {\phi}_i(x_i) \right) \det(g(s_j - x_i))_{1 \leq i,j \leq n} dx_1 \dots dx_n.
\end{equation*}

We now use Proposition \ref{proposition2} to write:
\begin{equation*}
    D = \int_{x_1 < \ldots < x_n} \underset{\sigma \in \mathfrak{S}_n}{\sum} \left( \overset{n}{\underset{i=1}{\prod}} f(x_{\sigma(i)}) \right) \left( \overset{n}{\underset{i=1}{\prod}} {\phi}_i (x_{\sigma(i)}) \right) \det (g(s_j - x_{\sigma (i) })) dx_1 \ldots dx_n.
\end{equation*}

For any permutation $\sigma \in \mathfrak{S}_n$, we have
$$\overset{n}{\underset{i=1}{\prod}} f(x_{\sigma(i)}) = \overset{n}{\underset{i=1}{\prod}} f(x_i).$$
    
Moreover, the function $g$ being in $PF_n$, we know that the function $d$ defined by
\begin{equation*}
    d(x_1, \ldots ,x_n) := \det (g(s_j - x_i))_{1 \leq i,j \leq n}.
\end{equation*}
is non-negative when $x_1 < \ldots < x_n$. Moreover, for any permutation $\sigma \in \mathfrak{S}_n$, we have
\begin{equation*}
    d(x_{\sigma(1)}, \ldots ,x_{\sigma(n)}) := \varepsilon(\sigma) \det (x_1, \ldots ,x_n).
\end{equation*}

We thus have:
\begin{align*}
    D & = \int_{x_1 < \ldots < x_n} \left( \overset{n}{\underset{i=1}{\prod}} f(x_i) \right) d(x_1, \ldots ,x_n) \left( \underset{\sigma \in \mathfrak{S}_n}{\sum} \varepsilon (\sigma) \overset{n}{\underset{i=1}{\prod}}  {\phi}_i(x_{\sigma(i)}) \right) dx_1 \ldots dx_n \\
    &= \int_{x_1 < \ldots < x_n} \left( \overset{n}{\underset{i=1}{\prod}} f(x_i) \right) d(x_1, \ldots ,x_n) \det({\phi}_i(x_j))_{1 \leq i,j \leq n} dx_1 \ldots dx_n \\
    &\geq 0.
\end{align*}
the inequality coming from $({\phi}_1, \ldots ,{\phi}_n)$ being in $GM_n$.
\end{proof}

An interesting consequence of Theorem \ref{theorem1} is the stability of the $PF_n$ class by convolution:
\begin{prop}
Let $X,Y$ be two independent random variables in the class $PF_n$. Then the variable $X+Y$ is also in $PF_n$.
\end{prop}
This result is already known, a proof can be found in \cite{karlin} or in \cite{schoenberg}. In the following proof, we suppose that at least one density is positive. 
\begin{proof}
Let $X,Y$ be two independent random variables in the class $PF_n$. We denote by $f$ and $g$ their densities and by $r = f\star g$ the density of the random variable $X+Y$. We assume that $f >0$. \\

We want to prove that for any $a_1 \leq \cdots \leq a_n$ and $b_1 \leq \cdots \leq b_n, ~ \det(r(a_i-b_j))_{1 \leq i,j \leq n}  \geq 0.$ \\
Let us fix an ordered n-uple $b_1 \leq \cdots \leq b_n$.
Now, we denote the function ${\phi}_i$ by
\begin{equation*}
{\phi}_i(x) := \dfrac{f(x-b_i)}{f(x)} .
\end{equation*}
We have that $f \in PF_n$, then by definition of the class $PF_n$ for any $a_1 \leq \ldots \leq a_n$:
$$\det({\phi}_i(a_j))_{1 \leq i,j \leq n} =  \frac{1}{\overset{n}{\underset{j=1}{\prod}} f(a_j)}  \det(f(a_j-b_i))_{1 \leq i,j \leq n} \geq 0.$$
In other words, $({\phi}_i)_i$ is in the class $GM_n$. So we can apply Theorem \ref{theorem1} to the $n$-uple of functions $({\phi}_i)_i$. We can also remark that:
$${\Phi}_i(a_j) = \dfrac{\int f(a_j-x) {\phi}_i(a_j-x) g(x) dx}{r(a_j)} = \dfrac{\int f(a_j - b_i -x) g(x) dx}{r(a_j)} = \dfrac{r(a_j-b_i)}{r(a_j)}.$$
Then by Theorem \ref{theorem1}, we have, for any $a_1 \leq \ldots \leq a_n$ :
$$ \det({\Phi}_i(a_j))_{1 \leq i,j \leq n} = \underset{j}{\prod} \dfrac{1}{r(a_j)}  \det(r(a_j-b_i))_{1 \leq i,j \leq n}  \geq 0.$$ 
Finally, we have proved that  for any
$a_1 \leq \ldots \leq a_n$ and $b_1 \leq \cdots \leq b_n$:
    $$\det(r(a_i-b_j))_{1 \leq i,j \leq n}  \geq 0.$$
\end{proof}

\section{Efron's Theorem with a stronger monotonicity assumption}
  We turn now to a refined version of Efron's theorem: under a stronger assumption on the monotonicity of the function $\phi$ on each variable we have a stronger result. We state it in two different frameworks, one for real random variables and a second for integer-valued random variables.
\subsection{Statement and examples}
The continuous version of the generalized Efron's theorem is stated as follows:
\begin{theo}
\label{theorem2}
Let X,Y be two independent log-concave real-valued random variables having densities. Let $a \ge 0$ be a real parameter and $\phi : \mathbb{R}^2 \rightarrow \mathbb{R}$ be a measurable function such that:
\begin{enumerate}
\item For any $ y \in \mathbb{R}$, the function $x \mapsto \exp(-ax) \phi(x,y)$ is non-decreasing.
\item For any $ x \in \mathbb{R}$, the function $y \mapsto \exp(-ay) \phi(x,y)$ is non-decreasing.
\end{enumerate}
 We define the function $\Phi : \mathbb{R} \rightarrow \mathbb{R}$ by
\begin{equation} \label{eq:defPhi}
\Phi(s) := \mathbb{E}[\phi(X,Y)|X+Y=s].
\end{equation}
We assume that $\Phi$ is well-defined. Then the function $s \mapsto \exp(-as) \Phi(s) $ is non-decreasing.
\end{theo}
If we take $a=0$ then conditions of Theorem \ref{theorem2} are equivalent to $\phi$ being a  non-decreasing function in each variable and Theorem \ref{theorem2} is equivalent to Efron's theorem. 
\medskip

Under additional regularity assumptions on $\phi$, Theorem~\ref{theorem2} can be restated as follows:
\begin{cor}
\label{corollary}
Let X,Y be two independent log-concave real random variables with densities. \\ 
Let $a \geq 0$ be a real parameter and $\phi : \mathbb{R}^2 \rightarrow \mathbb{R}$ be a function such that:
\begin{enumerate}
\item $\forall (x,y) \in \mathbb{R}^2 \ , \ {\partial}_1 \phi(x,y) \geq a \phi(x,y)$.
\item $\forall (x,y) \in \mathbb{R}^2 \ , \ {\partial}_2 \phi(x,y) \geq a \phi(x,y)$.
\end{enumerate}
If the function $\Phi : \mathbb{R} \rightarrow \mathbb{R}$ defined as in equation~\eqref{eq:defPhi} is differentiable, then it satisfies:
\begin{equation*}\forall s \in \mathbb{R}, \Phi'(s) \geq a \Phi(s).
\end{equation*}
\end{cor}

We give another equivalent statement for Theorem~\ref{theorem2}: if ${\phi}_1 : \mathbb{R}^2 \rightarrow \mathbb{R}$ is defined by 
$${\phi}_1(x,y) := \exp(a(x+y))$$
and if we set $\phi_2 = \phi$, then the conclusion of Theorem \ref{theorem2} is equivalent to saying that $({\Phi}_1,{\Phi}_2) \in GM_2$ where $\forall s \in \mathbb{R}, {\Phi}_i(s) = \mathbb{E}[{\phi}_i(X,Y)|X+Y=s]$.

\bigskip

We turn now to the statement of the generalized Efron's Theorem in a discrete setting:
\begin{theo}
\label{theorem2_2} Let $X,Y$ be two independent log-concave random variables on $\mathbb{Z}$. Let $a \geq 0$ be a real parameter and $\phi : \mathbb{Z}^2 \rightarrow \mathbb{\mathbb{R}}$ be a function such that:
\begin{enumerate}
\item For any $ k_2 \in \mathbb{Z}$, the function $k_1 \mapsto \exp(-ak_1) \phi(k_1,k_2)$ is non-decreasing.
\item For any $ k_1 \in \mathbb{Z}$, the function $k_2 \mapsto \exp(-a k_2) \phi(k_1,k_2)$ is non-decreasing.
\end{enumerate}

We define the function $\Phi : \mathbb{Z} \rightarrow \mathbb{R} $ by
\begin{equation*}
\Phi(s) := \mathbb{E}[\phi(X,Y)|X+Y=s].
\end{equation*}
We assume that $\Phi$ is well-defined. Then the function $s \mapsto \exp(-as) \Phi(s)$ is non-decreasing.
\end{theo}

We might notice that conditions on $\phi$ are equivalent to:
\begin{enumerate}
\item $\forall (x,y) \in \mathbb{Z}^2, ~ \phi(x+1,y) \geq e^a \phi(x,y)$,
\item $\forall (x,y) \in \mathbb{Z}^2, ~ \phi(x,y+1) \geq e^a \phi(x,y)$.
\end{enumerate}

Furthermore, if we denote ${\Delta}_1 \phi(x,y) = \phi(x+1,y) - \phi(x,y)$ and ${\Delta}_2 \phi(x,y) = \phi(x,y+1) - \phi(x,y)$ the discrete derivative operators, then the same conditions are also equivalent to:
\begin{enumerate}
\item $\forall (x,y) \in \mathbb{Z}^2, ~ {\Delta}_1 \phi(x,y) \geq (e^a-1) \phi(x,y)$,
\item $\forall (x,y) \in \mathbb{Z}^2, ~ {\Delta}_2 \phi(x,y) \geq (e^a-1) \phi(x,y)$.
\end{enumerate}
We can remark that, under this form, these conditions can be seen as a discrete version of conditions of Corollary \ref{corollary}.

\subsubsection{Examples}
We will consider, in this paragraph, examples for which variables are separated, i.e. functions $\phi : \mathbb{R}^2 \rightarrow \mathbb{R}$ that can be writtten as a product of two non-negative real functions:
\begin{equation*}
\forall (x,y) \in \mathbb{R}^2, \phi(x,y) = f(x) g(y).
\end{equation*}
We assume that:
\begin{enumerate}
    \item For any $y \in \mathbb{R}$, the function $x \mapsto \exp(-ax) f(x)$ is non-decreasing.
    \item For any $x \in \mathbb{R}$, the function $y \mapsto \exp(-ay) g(y)$ is non-decreasing.
\end{enumerate}
We denote the function $\Phi$ by
$$\Phi(s) := \mathbb{E}[\phi(X,Y) | X+Y=s].$$
We assume that the function $\Phi$ is well-defined. Then by Theorem \ref{theorem2}, we have that the function $s \mapsto \exp(-as) \Phi(s)$ is non-decreasing. 

\medskip

If $f$ and $g$ can be written under the form $f(x) = \exp(u(x))$ and $g(y) = \exp(v(y))$ where $u$ and $v$ are two real functions, then \begin{equation*}
    \forall (x,y) \in \mathbb{R}^2 \ , \ \phi(x,y)= \exp(u(x))  \exp(v(y)) = \exp(u(x) + v(y)).
\end{equation*} Then, conditions of Theorem \ref{theorem2} are equivalent to:
\begin{enumerate}
    \item For any $y \in \mathbb{R}$, the function $x \mapsto u(x) - ax$ is non-decreasing.
    \item For any $x \in \mathbb{R}$, the function $y \mapsto v(y) - ay$ is non-decreasing.
\end{enumerate}
For a non-trivial example, we can consider 
\begin{equation*} \phi(x,y) = \exp[x(x^2+ \alpha) + y(y^2 +\beta)].
\end{equation*}
By Theorem \ref{theorem2}, we have that the function $s \mapsto \exp(-\gamma s) \Phi(s)$ is  non-decreasing  with the choice of parameter $\gamma = \min(\alpha,\beta)$.

\subsection{Proof of Theorem 3 and 4}
Theorem \ref{theorem2} and Theorem \ref{theorem2_2} is a direct consequence of Efron's theorem. \\
We have assumed that $\Phi$ is well-defined. Then, $\mathbb{E}[\exp(-a(X+Y))\phi(X,Y)|X+Y=s]$ is also well-defined because we have that for any $s \in \mathbb{R}$:
\begin{equation*}
    \mathbb{E}[\exp(-a(X+Y))\phi(X,Y)|X+Y=s] = \mathbb{E}[\exp(-as)\phi(X,Y)|X+Y=s] = \exp(-as) \Phi(s).
\end{equation*}
We can remark that for any $(x,y) \in \mathbb{R}^2$ :
\begin{equation*}
 \exp(-a(x+y)) \phi(x,y) = \exp(-ax) \exp(-ay) \phi(x,y). 
\end{equation*}
 Then, the function $(x,y) \mapsto \exp(-a(x+y)) \phi(x,y)$ is non-decreasing in each variable because $x \mapsto \exp(-ax)\phi(x,y)$ and $y \mapsto \exp(-ay)\phi(x,y)$ are non-decreasing. \\
By the Efron's theorem, apply to the function $(x,y) \mapsto \exp(-a(x+y))\phi(x,y)$, we have that the function $s \mapsto \mathbb{E}[\exp(-a(X+Y))\phi(X,Y)|X+Y=s]$ is non-decreasing. But we have that :
\begin{equation*}
\mathbb{E}[\exp(-a(X+Y))\phi(X,Y)|X+Y=s] = \alpha(s) \Phi(s).
\end{equation*}
By definition of $\Phi$, we obtain $\mathbb{E}[\exp(-a(X+Y))\phi(X,Y)|X+Y=s] = \exp(-as) \Phi(s)$. \\
Finally, the function $s \mapsto \exp(-as) \Phi(s)$ is a non-decreasing function.

\section{Some generalizations of Theorem \ref{theorem2} and Theorem \ref{theorem2_2}}
In this part, we will show a more general result. We will adapt the result of the previous part for a more general real function $\alpha$. After, we will apply this result to some example. 
\begin{prop}
\label{efron_general}
Let $X,Y$ be two independent log-concave real-valued random variables. Let $\phi : \mathbb{R}^2 \rightarrow \mathbb{R}$ a measurable function and $\alpha : \mathbb{R} \rightarrow \mathbb{R}$ a measurable function such as :
\begin{enumerate}
\item For any $y \in \mathbb{R}$, the function $x \mapsto \alpha(x+y) \phi(x,y)$ is non-decreasing.
\item For any $x \in \mathbb{R}$, the function $y \mapsto \alpha(x+y) \phi(x,y)$ is non-decreasing.
\end{enumerate}
We define the function $\Phi : \mathbb{R} \rightarrow \mathbb{R}$ by :
\begin{equation*}
\Phi(s):= \mathbb{E}[\phi(X,Y)|X+Y=s].
\end{equation*}
We assume that $\Phi$ is well-defined. Then, we have that the function $s \mapsto \alpha(s) \Phi(s)$ is non-decreasing. \\ 
\end{prop}

\begin{proof}
By assumption, we have that $\Phi$ is well-defined. Then, we have that $\mathbb{E}[\alpha(X+Y) \phi(X,Y) | X+Y =s]$ is well-defined because :
$$\mathbb{E}[\alpha(X+Y) \phi(X,Y) | X+Y =s] = \mathbb{E}[\alpha(s) \phi(X,Y) | X+Y =s] = \alpha(s) \Phi(s).$$
We apply the Efron's theorem to the function $(x,y) \mapsto \alpha(x+y) \phi(x,y)$. Then we have that $s \mapsto \mathbb{E}[\alpha(X+Y) \phi(X,Y)|X+Y=s]$ is non-decreasing. But, we can remark that :
\begin{equation*}
 \forall s \in \mathbb{R}, \mathbb{E}[\alpha(X+Y) \phi(X,Y)|X+Y=s]= \mathbb{E}[\alpha(s) \phi(X,Y)|X+Y=s].
\end{equation*}
Finally, the function $s \mapsto \alpha(s) \Phi(s)$ is non-decreasing.
\end{proof}

\begin{cor}
Let $X,Y$ be two independent log-concave real-valued random variables. Let $f: \mathbb{R} \rightarrow \mathbb{R}$ and $g: \mathbb{R} \rightarrow \mathbb{R}$ be two real functions. Let $\alpha : \mathbb{R} \rightarrow \mathbb{R}$ a measurable function such as :
\begin{enumerate}
\item For any $\lambda \in \mathbb{R}$, the function $x \mapsto \alpha(x+\lambda) f(x)$ is non-decreasing.
\item For any $\lambda \in \mathbb{R}$, the function $y \mapsto \alpha(y + \lambda) g(y)$ is non-decreasing.
\end{enumerate}
We define the function $\Phi : \mathbb{R} \rightarrow \mathbb{R}$ by :
\begin{equation*}
\Phi(s):= \mathbb{E}[f(X)g(Y)|X+Y=s].
\end{equation*}
We assume that $\Phi$ is well-defined. Then, we have that the function $s \mapsto \alpha(s) \Phi(s)$ is non-decreasing. 
\end{cor}

\begin{proof}
To prove the corollary, we have only to apply the Proposition \ref{efron_general} to the function $\phi : (x,y) \mapsto f(x) g(y)$.
\end{proof}

We can remark that the Theorem 3 and 4 are a particular case of the proposition \ref{efron_general}. In the case of Theorem 3 and 4, we have that $\alpha(s) := \exp(-as)$. \\
The Efron's Theorem is the particular case where $\alpha(s):= 1$. \\
\medskip

Now, we will show an example of application of the Proposition \ref{efron_general}. \\
In this proposition, we will apply the proposition \ref{efron_general} to the function $\phi(x,y) := xy$ and $\alpha(s) := \dfrac{1}{s}$.
\begin{prop}
Let $X,Y$ be two independent log-concave real-valued random variables. We assume that $X,Y$ are positive.
We define the function $\Phi : \mathbb{R}_{+}^{*} \rightarrow \mathbb{R}$ by :
\begin{equation*}
\Phi(s):= \mathbb{E}[XY|X+Y=s].
\end{equation*}
We assume that $\Phi$ is well-defined. Then, we have that the function $s \mapsto \dfrac{\Phi(s)}{s}$ is non-decreasing on $\mathbb{R}_{+}^{*}$. 
\end{prop}

\begin{proof}
We define the function $\Psi : (\mathbb{R}_{+}^{*})^2 \rightarrow \mathbb{R}$ by :
\begin{equation*}
\Psi(x,y) := \dfrac{xy}{x+y}.
\end{equation*}
The function $\Psi$ is differentiable and we have that :
\begin{equation*}
\forall (x,y) \in (\mathbb{R}_{+}^{*})^2 , \dfrac{\partial \Psi}{\partial x}(x,y) = \dfrac{y^2}{(x+y)^2} \geq 0.
\end{equation*}

\begin{equation*}
\forall (x,y) \in (\mathbb{R}_{+}^{*})^2, \dfrac{\partial \Psi}{\partial y}(x,y) = \dfrac{x^2}{(x+y)^2} \geq 0.
\end{equation*}
Then, the function $\Psi$ is non-decreasing in each variable. By the proposition \ref{efron_general}, we obtain that the function $s \mapsto \dfrac{\Phi(s)}{s}$ is a non-decreasing function. 
\end{proof}

\end{document}